\documentclass{amsart}
\usepackage{amsmath, amsthm, amssymb}
\usepackage[all]{xy}

\theoremstyle{plain}
\newtheorem{theorem}{Theorem}[section]
\newtheorem{corollary}[theorem]{Corollary}
\newtheorem{lemma}[theorem]{Lemma}
\newtheorem{proposition}[theorem]{Proposition}

\theoremstyle{definition}
\newtheorem{definition}[theorem]{Definition}
\newtheorem{remark}[theorem]{Remark}

\theoremstyle{remark}

\begin{document}
\title[covering groupoids]
      {covering groupoids}
\author{Zhi-Ming Luo}
\address{Department of Mathematics\\
         University of Haifa\\
         Haifa, 31905, Israel}
\email{luo@math.haifa.ac.il}
\date{December 5, 2004}

\begin{abstract}
Topos properties of the category of covering groupoids over a fixed groupoid are discussed.
A classification result for connected covering groupoids over a fixed groupoid
analogous to the fundamental theorem of Galois theory is given.
 \par
\smallskip
\end{abstract}

\maketitle
%
%
%
\section{Introduction}
The theory of covering spaces
is one of the most beautiful theories in classic algebraic topology~\cite {R}.
After studying the relations between the fundamental groupoids of covering spaces and those of
base spaces, it is not hard to extend to the case of general groupoids and to give the
definition of covering groupoids. Since groupoid is an algebraic object and its $n$th homotopy
groups vanish for $n \geq 2$, it is relatively
easier than the case of spaces. Many mathematician deal with the theory of covering spaces
using that of covering groupoids (c.f. \cite {B}).

Brown \cite {B} and Higgins \cite {H} have elegantly treated with covering groupoids. They give
the definition of covering projection of groupoids, concisely prove covering homotopy theorem,
unique lifting theorem and existence theorem of covering
groupoids, gracefully apply the theory of
covering groupoids to group theory.

In this paper we regard a groupoid with the sets of all
arrows having same codomain on each object as a
Grothendieck site, the category of covering projections over such a groupoid is equivalent to
the category of sheaves on such a Grothendieck site. From the viewpoint of topos theory,
covering groupoids have some analogous properties of $\acute{\mbox{e}}$tale spaces rather than
those of covering spaces. We explicitly describe the subobject
classifier and exponentials in the category of covering projections in Sec.~2.

In Sec.~3 we study covering transformations and orbit morphisms of connected covering
projections. It is shown that connected covering groupoid over a base groupoid is equivalent to
an orbit groupoid which is a universal covering groupoid over such a base groupoid modulo some
covering transformation group. The family of subgroups of covering transformation group of a
universal covering projection and the family of equivalence classes of connected covering
groupoids form two lattices. We have a classification result (Theorem~\ref {T: main}) for the
equivalence classes of connected covering groupoids which say that there is an
order-reversing bijection between the two lattices.


\section{Definitions and basic properties}

A \emph{groupoid} $\mathbf{G}$ is a category in which every arrow is invertible.
Given an object $G$ in the groupoid $ \mathbf{G}$, a \emph{sieve} on $G$ is a set
$S$ of arrows with codomain $G$ such that
\begin{center}
$f \in S$ and the composite $fh$ is defined implies $fh \in S$.
\end{center}
Since every arrow in $ \mathbf{G}$ is invertible, then every nonempty sieve on $G$
is the set of all arrows into $G$, i.e., the \emph{maximal sieve} $t(G)$ on $G$.

We call the maximal sieve $t(G)$ as a \emph{neighborhood} of $G$ in the groupoid $ \mathbf{G}$.

\begin{definition}
Let $p: \mathbf{\tilde{G}} \to \mathbf{G}$ be a morphism of groupoids. An ordered pair
$( \mathbf{\tilde{G}}, p)$ is a \emph{covering groupoid} if for each object $ \tilde{G}$
of $ \mathbf{\tilde{G}}$ the restriction of $p$
$$t( \tilde{G}) \to t(p\tilde{G})$$
is bijection. The morphism $p$ is called the \emph{covering projection}. The covering
projection $p$ is called \emph{connected} if both $ \mathbf{\tilde{G}}$ and
$ \mathbf{G}$ are connected.
\end{definition}

\begin{remark}
The maximal sieve $t(G)$ is same as the co-star in the sense of Higgins~\cite {H} and
equivalent
to the star $St_{ \mathbf{G}} G$ in the sense of Brown~\cite {B} and Higgins~\cite {H}
in the case of groupoids. Hence their definition coincides with ours.
\end{remark}

Since $p: t( \tilde{G}) \to t(p\tilde{G})$ is bijection, for every $a \in t(p\tilde{G})$
there is a unique
$ \tilde{a} \in t( \tilde{G})$ with $p( \tilde{a}) = a$. We call $ \tilde{a}$ as a
\emph{lifting} of $a$. When $p$ is a connected covering projection it is surjective on objects.

Let $ \mathbf{G}$ be a groupoid, the category $ \mathcal{COV}( \mathbf{G})$ of covering
projections of $ \mathbf{G}$ has as objects the covering projections
 $p: \mathbf{H} \to \mathbf{G}$
and has arrows the commutative diagrams of morphisms of groupoids,
$$\xymatrix{
\mathbf{H}\ar[rr]^-f\ar[rd]_-p&& \mathbf{K}\ar[dl]^-q\\
&\mathbf{G}}
$$
where $p$ and $q$ are covering projections. It is convenient to just write such an arrow as $f$.
The groupoid-morphism $f: \mathbf{H} \to \mathbf{K}$ is also a covering projection since for
every object $H \in \mbox{Ob}( \mathbf{H})$ there exist bijections $t(H) \cong t(pH) \cong
t(qfH) \cong t(fH)$.

For any groupoid $ \mathbf{G}$, denote the set of connected components of $ \mathbf{G}$ by
$ \pi_0(\mathbf{G})$. For any object $G \in \mbox{Ob}(\mathbf{G})$, all arrows $G \to G$ in
$ \mathbf{G}$ form a group. We denote this group by $ \pi(\mathbf{G}, G)$ and call it as the
\emph{fundamental group} of $ \mathbf{G}$ at $G$ (Higgins calls it as vertex group at $G$ and
Brown calls it as object group at $G$). Let $p: \mathbf{\tilde{G}}\to \mathbf{G}$ be a covering
projection, it is clear that $\pi( \mathbf{\tilde{G}}, \tilde{G})$ is isomorphic to the
subgroup
$p_\ast\pi( \mathbf{\tilde{G}}, \tilde{G})$  of $\pi( \mathbf{{G}}, p\tilde{G})$.

Given a groupoid-morphism $f: \mathbf{F} \to \mathbf{G}$, if $f$ induces a bijection $f_\ast:
\pi_0( \mathbf{F}) \to \pi_0( \mathbf{G})$ and an isomorphism $f_\ast: \pi( \mathbf{F}, F)
\to \pi( \mathbf{G}, fF)$ for each $F \in \mbox{Ob}( \mathbf{F})$, we call $f$ as a \emph{
weak equivalence}. A groupoid-morphism is a weak equivalence and a covering projection if and
only if it is an isomorphism.

Higgins (Corollary 35 \cite{H}) and Brown (Theorem 9.3.3 \cite{B}) prove the following unique
lifting theorem:

\begin{theorem}[Unique lifting theorem]
Let $( \mathbf{\tilde{G}}, p)$ be a covering groupoid of $ \mathbf{G}$ with $p\tilde{G} = G$
where $\tilde{G} \in \mbox{Ob}( \mathbf{\tilde{G}}$) and $G \in \mbox{Ob}( \mathbf{G})$, and
$f: \mathbf{F} \to \mathbf{G}$ a groupoid morphism with $f(F) = G$ where $F \in \mbox{Ob}(F)$
such that $ \mathbf{F}$ is connected. Then $f$ lifts to a morphism $ \tilde{f}: \mathbf{F} \to
\mathbf{\tilde{G}}$ with $ \tilde{f}(F) = \tilde{G}$ if and only if $f_\ast\pi( \mathbf{F}, F)
\subset p_\ast\pi( \mathbf{\tilde{G}}, \tilde{G})$, and if this lifting exists, then it is
unique.
\end{theorem}

For every subgroup $\Gamma$ of the fundamental group $\pi( \mathbf{G}, G)$
of $ \mathbf{G}$ at $G$,
Higgins \cite{H} and Brown (Theorem 9.4.3 \cite{B}) also prove the following covering-groupoid
existence theorem:

\begin{theorem}[Existence]
Let $G$ be an object of the connected groupoid $ \mathbf{G}$, and let $\Gamma$ be a subgroup
of the fundamental group $\pi( \mathbf{G}, G)$, then there is a covering groupoid
$( \mathbf{\tilde{G}}, p)$ of $ \mathbf{G}$ with $p( \tilde{G}) = G$ where $ \tilde{G}
\in \mbox{Ob}( \mathbf{\tilde{G}}) $ and $p_\ast\pi( \mathbf{\tilde{G}}, \tilde{G}) = \Gamma$.
\end{theorem}

Let $ \mathbf{G}$ be a connected groupoid, there exists a connected covering groupoid
$ \mathbf{\tilde{G}}$ of $ \mathbf{G}$ with trivial fundamental group from the existence
theorem. For every covering groupoid $ \mathbf{H}$ of $ \mathbf{G}$,
$ \mathbf{\tilde{G}}$ is a covering groupoid of $ \mathbf{H}$ from the unique lifting theorem.
This groupoid $ \mathbf{\tilde{G}}$ is called a \emph{universal covering groupoid} of
$ \mathbf{G}$.

Let $p: \mathbf{\tilde{G}} \to \mathbf{G}$ be a covering projection, each component of
$ \mathbf{\tilde{G}}$ covers a components of $ \mathbf{G}$. For every component of
$ \mathbf{G}$ its preimage is the union of components in $ \mathbf{\tilde{G}}$ covering it.
Without lost of generality we assume that the base groupoid $\mathbf{G}$ is connected.

Write $ \mathcal{CONCOV}( \mathbf{G})$ for the subcategory of $ \mathcal{COV}( \mathbf{G})$
in which the objects are connected covering projections.

\smallskip

Given a Grothendieck topology $J$ on the groupoid $ \mathbf{G}$ defined by
\begin{center}
$S \in J(G)$ iff the sieve $S$ is nonempty.
\end{center}
i.e., $J$ is the atomic topology on $ \mathbf{G}$ \cite {M-M}, $J$ is the minimal
topology \cite {J} (there are only two topologies on the groupoid $ \mathbf{G}$ --- the maximal
topology and the minimal topology). Every presheaf on the site $( \mathbf{G}, J)$
is a sheaf, i.e., $\bf {Sets}^{{\mathbf{G}}^{op}} = {\bf Sh}( \mathbf{G}, J)$.

Each covering projection $p: \mathbf{\tilde{G}} \to \mathbf{G}$ over $ \mathbf{G}$
determines a $ \mathbf{G}$-indexed family $\{ \mathbf{H}_G | G \in \mathbf{G} \}$
of groupoids, consisting of the faithful and full subgroupoids of $ \mathbf{\tilde{G}}$
$$ \mathbf{H}_G = p^{-1}\{G\},$$
its objects are the elements in the set $p_{\mbox{ob}}^{-1}\{G\}$ in $ \mathbf{\tilde{G}}$,
where $p_{\mbox{ob}}$ is the restriction of $p$ on objects
$p_{\mbox{ob}}: \mbox{ob}( \mathbf{\tilde{G}} ) \to \mbox{ob}( \mathbf{G})$.

For each arrow $f: D \to C$ in $ \mathbf{G}$ and each object $ \tilde{C} \in \mathbf{H}_C$,
the map $t( \tilde{C}) \to t(C)$ is bijective, hence there exists a unique arrow
$ f_{ \tilde{C}}: \tilde{C} \to \tilde{D}$ in $ \mathbf{\tilde{G}}$ such that $ \tilde{D}
\in \mathbf{H}_D$ and $ (pf_{ \tilde{C}})^{-1} = f$.
For each arrow $s: \tilde{C} \to \tilde{C'}$
in $ \mathbf{H}_C$ there exists a unique commutative
diagram
$$\xymatrix{
\tilde{C}\ar[r]^s\ar[d]_{f_{\tilde{C}}} & \tilde{C^\prime}\ar[d]^{f_{\tilde{C^\prime}}}\\
\tilde{D}\ar[r]_t & \tilde{D^\prime}. }$$
Define a map $f_s: s \to t$. It's easy to verify that $ \mathbf{H}_f = \{f_{\tilde{C}},
f_s\}$ is a morphism $ \mathbf{H}_f: \mathbf{H}_C \to \mathbf{H}_D$ of groupoids.
$f$ is invertible, so is $ \mathbf{H}_f$, hence groupoids $ \mathbf{H}_C$ and $ \mathbf{H}_D$
are isomorphic.

Given a covering groupoid $ (\mathbf{\tilde{G}}, p)$, a contravariant
functor $ \mathbf{H}$ is determined
from the groupoid $ \mathbf{G}$ to the category of groupoids. The action of $ \mathbf{H}$
on the object $G$ of $ \mathbf{G}$ is the groupoid $ \mathbf{H}_G$ and on the arrow $f$ of
$ \mathbf{G}$ is the groupoid-morphism $ \mathbf{H}_f$. Hence $ (\mathbf{\tilde{G}}, p)$
determines a sheaf of groupoids $ \mathbf{H}$ on the site $( \mathbf{G}, J)$.

\begin{proposition}\label {P: sheaf}
Each covering groupoid of the groupoid $ \mathbf{G}$ corresponds a sheaf of groupoids on
the site $( \mathbf{G}, J)$.
\end{proposition}

Let $p: \mathbf{\tilde{G}} \to \mathbf{G}$ be a connected covering projection, the groupoids
$ \mathbf{H}_C$ and $ \mathbf{H}_D$ are isomorphic ( call them as fibres over $C$ and $D$,
respectively). Hence there is a bijection between the sets $\mbox{Ob}(\mathbf{H}_C)$  and
$\mbox{Ob}(\mathbf{H}_D)$, they have same cardinal denoted by $|\mbox{Ob}(\mathbf{H}_C)|$.
When the cardinal is a finite number $n$, we call $p$ an \emph{$n$-fold} covering.

\begin{remark}
Higgins (Proposition 30 \cite{H}) shows that the category $ \mathcal{COV}( \mathbf{G})$ is
equivalent to the functor category $\bf{Sets}^{ \mathbf{G}}$. Let $ \mathbf{G}^{op}$
be the opposite of $ \mathbf{G}$, the category $\mathcal{COV}( \mathbf{G}^{op})$ is
equivalent to the category $\mathcal{COV}(\mathbf{G})$, hence  $\mathcal{COV}(\mathbf{G})$
is equivalent to $\bf{Sets}^{ \mathbf{G}^{op}}$ (it is $\bf Sh( \mathbf{G},
J)$). So $\mathcal{COV}(\mathbf{G})$ is a Grothendieck topos. From Proposition~\ref {P: sheaf}
it is clear that $\bf Sh( \mathbf{G}, J)$ is equivalent to a subcategory of the category of
sheaves of groupoids $\bf GpdSh( \mathbf{G}, J)$ on the site $( \mathbf{G}, J)$.
\end{remark}

For any topological space $X$, it is well-known that there is an equivalence of categories
$$\Lambda: {\bf Sh}(X) \rightleftarrows {\bf Etale}(X): \Gamma$$
where ${\bf Sh}(X)$ is the category of all sheaves over $X$, ${\bf Etale}(X)$ is the category
of all $\acute{\mbox{e}}$tale bundles over $X$, $\Lambda$ assigns to each sheaf $P$ the bundle
of germs of $P$ and $\Gamma$ assigns to each bundle $p: Y \to X$ the sheaf of all
cross-sections of $Y$. From above remark we know that there is an equivalence of categories
$$ {\bf Sh}( \mathbf{G}, J) \rightleftarrows \mathcal{COV}( \mathbf{G})$$
From this comparison category $\mathcal{COV}( \mathbf{G})$ is similar to category
${\bf Etale}(X)$. From the definition of $\acute{\mbox{e}}$tale map $p: E \to X$, $p$ is
a local homeomorphism: To each $e \in E$ there exists an open set $V$ with $e \in V \subset E$
such that $pV$ is open in $X$ and $p|_V$ is a homeomorphism $V \to pV$. But the definition
of groupoid covering projection $p: \mathbf{\tilde{G}} \to \mathbf{G}$ requires $p$ is also
``local" bijection: To each $ \tilde{G} \in \mbox{Ob}( \mathbf{\tilde{G}}) $ and its neighbor
$t( \tilde{G})$ there is a bijection $t( \tilde{G}) \to t( p\tilde{G})$. Hence covering
groupoids have some analogous properties of $\acute{\mbox{e}}$tale spaces of $X$ rather than
those of covering spaces of $X$.

The category $ \mathcal{COV}( \mathbf{G}) $ has initial object the empty covering $0: \emptyset
\to \mathbf{G}$ and terminal object the identity covering $1: \mathbf{G} \to \mathbf{G}$. Since
$ \mathcal{COV}( \mathbf{G}) $ is a topos, it is closed under all pushouts and pullbacks.
$ \mathcal{CONCOV}(\mathbf{G})$ is a subcategory of $ \mathcal{COV}( \mathbf{G})$, hence in
 $ \mathcal{CONCOV}(\mathbf{G})$ any two morphisms with common domain have
 pushout in $ \mathcal{COV}( \mathbf{G})$. Since these morphisms are connected
 covering projection, the pushouts is also connected, hence
 $ \mathcal{CONCOV}(\mathbf{G})$ is closed under pushouts. But it does not need to be closed
 under pullbacks since pullbacks may be non-connected. For example, let $\mathbf{G}$
 be a groupoid with nontrivial fundamental groups and $(\mathbf{\tilde{G}}, p)$ a connected
 universal covering groupoid, then the pullback of $p$ and itself is not connected. Hence
 $ \mathcal{CONCOV}(\mathbf{G})$ is not a topos.

The presheaf category ${\bf Sets}^{ \mathbf{G}^{op}}$ has a subobject classifier
$\Omega$ defined on objects by
$$
\begin{aligned}
\Omega(G)
&=\{S~~|~~S ~~\mbox{is a sieve on} ~~G ~~\mbox{in} ~~\mathbf{G}\}\\
&= \{t(G), ~~\emptyset\}
    \end{aligned}
$$
The universal monic ``true" morphism is the natural transformation
$$\mbox{true}: 1 \to \Omega$$
such a ``true" assigns to each object $G$ of $ \mathbf{G}$ a function $\mbox{true}_G: 1(G)
\to \Omega(G)$ with $\mbox{true}_G(\ast) = t(G)$ where $\ast$ is the unique element in $1(G)$.

Translate into the language of $ \mathcal{COV}( \mathbf{G})$, the category
$ \mathcal{COV}( \mathbf{G})$ has a subobject classifier the covering projection
$$\Omega = 1 \sqcup 1: \mathbf{G} \sqcup \mathbf{G} \to \mathbf{G} $$
where $\mathbf{G} \sqcup \mathbf{G}$ is the disjoint union of $\mathbf{G}$ and itself.

The universal monic ``true" morphism is the inclusion covering projection
$$\mbox{true}: \mathbf{G} \to \mathbf{G} \sqcup \mathbf{G}$$

Let $( \mathbf{H}, q)$ be a covering groupoid of $ \mathbf{G}$, for each monic $ \mathbf{S}
\rightarrowtail \mathbf{H}$ in $ \mathcal{COV}( \mathbf{G})$ there is a unique arrow $\phi$
which, with the monic $\mbox{true}: \mathbf{G} \to \mathbf{G} \sqcup \mathbf{G}$, forms
a pullback square
$$
\xymatrix{ \mathbf{S}\ar[r]\ar@{>->}[d]& \mathbf{G}\ar[d]^{\mbox{true}}\\
\mathbf{H}\ar@{.>}[r]_-\phi& \mathbf{G} \sqcup \mathbf{G}}
$$
hence all subobjects of $\mathbf{H}$ are the equivalence classes of the unions of its
connected
components. When $\mathbf{H}$ is connected, it has just two subobjects $\emptyset$ and the
equivalence class of itself.

If $p: \mathbf{\tilde{G}} \to \mathbf{G}$ be a nonempty connected covering projection,
then $p$
is epi. Since if there exist two groupoid-morphisms $s, t: \mathbf{G} \rightrightarrows
 \mathbf{K}$ with $sp = tp$, for every $G \in \mbox{Ob}( \mathbf{G})$ there is $\tilde{G} \in
 \mbox{Ob}( \mathbf{\tilde{G}})$ with $p( \tilde{G}) = G$ since $ \mathbf{G}$ is connected and
 $ \mathbf{\tilde{G}} $ is nonempty, hence $s(G) = s(p( \tilde{G})) = t(p( \tilde{G})) = t(G)$.
 For every morphism $a \in \mathbf{G}(G, G')$ it has a lifting $ \tilde{a} \in
 \mathbf{\tilde{G}}( \tilde{G}, \tilde{G}')$ with $p( \tilde{G}) = G, p( \tilde{G}') = G'$ and
 $p(\tilde{a} ) = a$, hence $s(a) = s(p( \tilde{a})) = t(p( \tilde{a})) = t(a)$, so $s = t$.

When $ \mathbf{S}$ is nonempty, $ \mathbf{H}$ is connected and $s: \mathbf{S} \rightarrowtail
\mathbf{H}$ is monic, $s$ is epi since $s$ is a nonempty connected covering projection, hence
$s$ is an isomorphism. This coincides with above argument. The category
$ \mathcal{CONCOV}(\mathbf{G})$ has not subobject classifier.

A topos $ \mathcal{E}$ is said to be \emph{Boolean} iff for every object $E$ of $ \mathcal{E}$
the partially ordered set $\mbox{Sub}(E)$ of subobjects of $E$ is a Boolean algebra.

For every covering groupoid $ ( \mathbf{H}, q)$ in the topos $ \mathcal{COV}( \mathbf{G})$,
the partially ordered set $\mbox{Sub}( \mathbf{H})$ is isomorphic to the power set
$\mbox{P}( \mathbf{H})$ of the set of connected components of $ \mathbf{H}$. Power set is
Boolean, so

\begin{proposition}
The topos $ \mathcal{COV}( \mathbf{G})$ is Boolean.
\end{proposition}

The topos ${\bf Sets}^{ \mathbf{G}^{op}}$ has exponentials. Assume that $Q, P, R \in
{\bf Sets}^{ \mathbf{G}^{op}}$ and ${\bf y}C = {\mbox{Hom}}_{ \mathbf{G}}(-, C)$ the
representable functor, the presheaf $Q^P$ exists by Proposition~I.6.1. \cite {M-M} such that
${\mbox{Hom}}_{ \mathbf{\hat{G}} }(R \times P, Q) \cong
{\mbox{Hom}}_{ \mathbf{\hat{G}} }(R, Q^P)$ where ${ \mathbf{\hat{G}} }$ is just the topos
${\bf Sets}^{ \mathbf{G}^{op}}$. The formula for $Q^P$ at $C$ is
$$
\begin{aligned}
Q^P(C)
&={\mbox{Hom}}_{ \mathbf{\hat{G}} }({\bf y}C, Q^P)\\
&= {\mbox{Hom}}_{ \mathbf{\hat{G}} }({\bf y}C \times P, Q)
    \end{aligned}
$$
i.e., $Q^P(C)$ is the set of all natural transformations $\theta:
{\mbox{Hom}}_{ \mathbf{G}}(-, C) \times P \to Q$.

Given an arrow $g: D \to C$ in $ \mathbf{G}$, since $g$ has an inverse $g^{-1}: C \to D$, the
following maps
$$P(g): P(C) \to P(D),$$
$$Q(g): Q(C) \to Q(D),$$
$${\mbox{Hom}}(g, C): {\mbox{Hom}}_{ \mathbf{G}}(C, C) \to {\mbox{Hom}}_{ \mathbf{G}}(D, C)$$
are bijections, where ${\mbox{Hom}}(g, C)(f) = fg$ when
$f \in {\mbox{Hom}}_{ \mathbf{G}}(C, C)$.

The natural transformations $\theta$ correspond exactly to the sets $\{ \theta_C, C \in
\mbox{Ob}( \mathbf{G})\}$ such that for every arrow $g: D \to C$ in $ \mathbf{G}$,
the following diagram is commutative
$$\xymatrix{
{\mbox{Hom}}_{ \mathbf{G}}(C, C) \times P(C)\ar[r]^-{\theta_C}
\ar[d]_-{\mbox{Hom}(g, C) \times P(g)}&Q(C)\ar[d]^-{Q(g)}\\
{\mbox{Hom}}_{ \mathbf{G}}(D, C) \times P(D)\ar[r]_-{\theta_D}&Q(D)}
$$

Given an element $(s,d) \in {\mbox{Hom}}_{ \mathbf{G}}(D, C) \times P(D)$, since
$\mbox{Hom}(g, C)$ and $ P(g)$ are bijective, $\theta_D(s,d) = Q(g)(\theta_C(sg^{-1},
(P(g))^{-1}d)$, but $\theta_D$ is a map independently on $g$, hence it is determined by
$\theta_C$. So $\theta$ corresponds exactly to $\theta_C$.

Given $g: C \to C$ in ${\mbox{Hom}}_{ \mathbf{G}}(C, C)$ where
${\mbox{Hom}}_{ \mathbf{G}}(C, C)$ is just the fundamental group of $ \mathbf{G}$ at $C$, the
sets of all $P(g)$ and $Q(g)$ can be regarded as the right action of
${\mbox{Hom}}_{ \mathbf{G}}(C, C)$ on $P(C)$ and $Q(C)$, respectively. Given $x \in P(C)$
and $y \in Q(C)$, write $P(g)x $ as $ xg$ and $Q(g)y$ as $yg$, respectively.

Similarly, there exists a commutative diagram
$$\xymatrix{
{\mbox{Hom}}_{ \mathbf{G}}(C, C) \times P(C)\ar[r]^-{\theta_C}
\ar[d]_-{\mbox{Hom}(g, C) \times P(g)}&Q(C)\ar[d]^-{Q(g)}\\
{\mbox{Hom}}_{ \mathbf{G}}(C, C) \times P(C)\ar[r]_-{\theta_C}&Q(C)}
$$
this diagram means that $\theta_C$ is an equivariant ${\mbox{Hom}}_{ \mathbf{G}}(C, C)$-map
on ${\mbox{Hom}}_{ \mathbf{G}}(C, C) \times P(C)$: take element $(f, x) \in
{\mbox{Hom}}_{ \mathbf{G}}(C, C) \times P(C)$, then $\theta_C(f, x)g = \theta_C(fg, xg)$
where ${\mbox{Hom}}_{ \mathbf{G}}(C, C)$ acts  right on itself by multiplication.

$Q^P$ is a contravariant functor, for $g: C \to C$, write $g^\ast$ for $Q^P(g): Q^P(C) \to
Q^P(C)$, then $g^\ast\theta$ is the composite
$$\xymatrix{
{\mbox{Hom}}_{ \mathbf{G}}(-, C) \times P \ar[r]^-\theta & Q\\
{\mbox{Hom}}_{ \mathbf{G}}(-, C) \times P \ar[u]^-{{\mbox{Hom}}(1, g) \times 1}
\ar[ur]_-{g^\ast\theta}}
$$

For the corresponding case of $\theta_C$, $g^\ast\theta_C$ is the composite
$$\xymatrix{
{\mbox{Hom}}_{ \mathbf{G}}(C, C) \times P(C) \ar[r]^-{\theta_C} & Q(C)\\
{\mbox{Hom}}_{ \mathbf{G}}(C, C) \times P(C) \ar[u]^-{{\mbox{Hom}}(1, g) \times 1}
\ar[ur]_-{g^\ast\theta_C}}
$$
$g^\ast\theta_C$ can be regarded as the right action of $g$ on $\theta_C$, write $\theta_cg$
for $g^\ast\theta_C$. Taken element $(f, x) \in
{\mbox{Hom}}_{ \mathbf{G}}(C, C) \times P(C)$, then $\theta_Cg(f, x) = g^\ast\theta_C(f, x)
= \theta_C(\mbox{Hom}(1, g)\times 1(f, x)) = \theta_C(gf, x)$.

Let $\alpha: P(C) \to Q(C)$ be defined by $\alpha(x) = \theta_C(1, x)$ for $x \in P(C)$, then
$\theta_C(f, x) = \theta_C(1, xf^{-1})f = \alpha(xf^{-1})f$ for $(f, x) \in
{\mbox{Hom}}_{ \mathbf{G}}(C, C) \times P(C)$.

Given any map $\alpha: P(C) \to Q(C)$, let $\theta_C:
{\mbox{Hom}}_{ \mathbf{G}}(C, C) \times P(C) \to Q(C)$ be defined by $
\theta_C(f, x) = \alpha(xf^{-1})f$, hence $\theta_C(fg, xg) = \alpha(xg(fg)^{-1})fg
= \alpha(xf^{-1})fg = \theta_C(f,x)g$, so $\theta_C$ satisfies the equivariant map property.

$\theta_C$ corresponds exactly to $\alpha$, hence
$$Q^P(C) = \{\alpha~|~{\mbox{all maps}} ~\alpha: P(C) \to Q(C)\}$$
The group ${\mbox{Hom}}_{ \mathbf{G}}(C, C)$ acts right on $Q^P(C)$: for $g \in
{\mbox{Hom}}_{ \mathbf{G}}(C, C)$ and $\alpha \in Q^P(C)$, $\alpha g$ is a map such that
for $x \in P(C)$
$$
\begin{aligned}
\alpha g(x)
&=\theta_Cg(1, x)
=\theta_C(g, x)
=\theta_C(1, xg^{-1})g \\
&=\alpha(xg^{-1})g
    \end{aligned}
$$

Given $g: D \to C$, there is a map $Q^P(g): Q^P(C) \to Q^P(D)$ with $Q^P(g)(\alpha)$ is the
bottom arrow in the following diagram
$$\xymatrix{
P(C)\ar[r]^-\alpha\ar[d]_-{P(g)}&Q(C)\ar[d]^-{Q(g)}\\
P(D)\ar[r]_-{Q^P(g)(\alpha)}&Q(D)}
$$
i.e., $Q^P(g)(\alpha) = Q(g)\alpha(P(g))^{-1}$ (both $P(g)$ and $Q(g)$ are bijective). When
$g \in {\mbox{Hom}}_{ \mathbf{G}}(C, C)$, $Q^P(g)(\alpha)$ is just $\alpha g$.

Now we return to the case of covering groupoids. Given two covering projections of groupoids
$p: \mathbf{H} \to \mathbf{G}$ and $q: \mathbf{K} \to \mathbf{G}$, the exponential $p^q:
\mathbf{H}^\mathbf{K} \to \mathbf{G}$ is defined as follows. The groupoid
$\mathbf{H}^\mathbf{K}$ has as the objects the maps $\alpha: \mbox{Ob}( \mathbf{K}_G)
\to \mbox{Ob}( \mathbf{H}_G)$ and as arrows $g_\ast: \alpha \to \mathbf{H}_g \alpha
(\mathbf{K}_g)^{-1}$. The covering projection is $p^q( \alpha) = G$ and $p^q(g_\ast) = g$.

If we take every $\alpha$ mapping to one point, it is easy to know that $ \mathbf{H}$ is a
subgroupoid of $ \mathbf{H}^\mathbf{K}$. The simplest case is $( \mathbf{G}^\mathbf{G}, 1^1)
=( \mathbf{G}, 1)$. Both $ \mathbf{H}$ and $\mathbf{K}$ are connected, but
$ \mathbf{H}^\mathbf{K}$ does not need to be connected, hence the category
$ \mathcal{CONCOV}(\mathbf{G})$ has not exponentials.

\begin{theorem}
Let $f: \mathbf{H} \to \mathbf{G}$ be a morphism of groupoids, the pullback $p'$ of the
covering projection $p$ along the morphism $f$
$$
\xymatrix{
\mathbf{\tilde{H}}\ar[r]^-{f^\ast}\ar[d]_-{p'}&\mathbf{\tilde{G}}\ar[d]^p\\
\mathbf{H}\ar[r]_-f&\mathbf{G}}
$$
is a covering projection. $f$ induces a morphism of topoi $f^\ast:  \mathcal{COV}(\mathbf{G})
\to \mathcal{COV}(\mathbf{H})$ such that the subcategory generated by
$f^\ast(\mathcal{COV}(\mathbf{G}))$ is
equivalent to the functor category generated by $Ff$ where $Ff$ is a functor
$\xymatrix@1{\mathbf{H}^{op}\ar[r]^-f&\mathbf{G}^{op}\ar[r]^-F&\mathbf{Sets}}$.
\end{theorem}
\begin{proof}
For every $\tilde{H} \in \mbox{Ob}(\mathbf{\tilde{H}})$, it is easy to see that $t(\tilde{H})
\cong t(p'(\tilde{H}))$ from the definition of pullback of groupoids. There is a commutative
diagram
$$\xymatrix{
 \mathcal{COV}(\mathbf{G}) \ar[r]^-{f^\ast}\ar[d]_-\simeq & \mathcal{COV}(\mathbf{H})
 \ar[d]^-\simeq\\
{\bf Sets}^{ \mathbf{G}^{op}}\ar[r]_-{f^\ast}&{\bf Sets}^{ \mathbf{H}^{op}}}
$$
up to isomorphism. Both vertical arrows are equivalent, so are the images of both horizontal
arrows, then the categories which they generate are equivalent.
\end{proof}


\section{Covering transformation}

It is well-known that the theory of covering spaces is analogous to Galois theory. The theory
of covering groupoids has same analogy. In this section we give the detailed description.
Suppose that $F$ is a subfield of a field $E$. Recall that
$$\mbox{Gal}(E/F) = \{~\mbox{ automorphisms }~\sigma: E \to E~|~\sigma~\mbox{fixes} ~F~
\mbox{pointwise}\}.$$
If $i: F \hookrightarrow E$ is the inclusion, then an automorphism $\sigma$ of $E$ lies in
$\mbox{Gal}(E/F)$ if and only if the following diagram commutes:
$$\xymatrix{
E \ar[rr]^-\sigma&& E.\\
&F\ar[ul]^-i\ar[ur]_-i}
$$
In the case of covering groupoids, we expect that there is an isomorphism
$h: \mathbf{\tilde{G}} \to \mathbf{\tilde{G}}$ which lies in a group if
and only if the following diagram commutes:
$$\xymatrix{
\mathbf{\tilde{G}}\ar[dr]_p \ar[rr]^-h&& \mathbf{\tilde{G}}\ar[dl]^-p.\\
&\mathbf{G}}
$$
and such a group has analogous property of $\mbox{Gal}(E/F)$. Since the arrows in second
diagram are reversed to those in first diagram, we can regard the theory of covering groupoids
as a ``co-Galois theory". In this analogy, universal covering groupoids will play the role of
algebraic closures. In this section we restrict all covering projections to be connected.

\begin{theorem}\label {T: action}
Let $( \mathbf{\tilde{G}}, p)$ be a connected covering groupoid of $ \mathbf{G}$, let $G \in
{\mbox{Ob}}( \mathbf{G})$, and let $ \mathbf{H}_G$ be the fibre over $G$. Then
\begin{itemize}
\item [(i)] $\pi( \mathbf{G}, G)$ acts transitively on $\mbox{Ob}( \mathbf{H}_G)$.
\item [(ii)] If $ \tilde{G} \in \mbox{Ob}( \mathbf{H}_G)$, then the stabilizer of $ \tilde{G}$
is $ p_\ast\pi( \mathbf{\tilde{G}}, \tilde{G})$.
\item [(iii)] $| \mbox{Ob}( \mathbf{H}_G) | = [ \pi( \mathbf{G}, G):
p_\ast\pi( \mathbf{\tilde{G}}, \tilde{G})]$.
\end{itemize}
\end{theorem}
\begin{proof} (i) Define a function $ \mbox{Ob}( \mathbf{H}_G) \times  \pi( \mathbf{G}, G)
\to \mbox{Ob}( \mathbf{H}_G)$, denoted by $( \tilde{G}, f) \mapsto \tilde{G}f $. Since there
is a unique arrow $ \tilde{f} \in t( \tilde{G})$ with $p\tilde{f} = f$, let $\tilde{G}f
= \mbox{dom}( \tilde{f})$ the domain of arrow $ \tilde{f}$. The identity arrow $ \tilde{e}
\in t( \tilde{G})$ corresponds to the identity arrow $ e \in \pi( \mathbf{G}, G)$, so
$\tilde{G}e = \tilde{G}$. For any two arrows $f, g \in \pi( \mathbf{G}, G)$, there exist
unique arrows $ \tilde{f} \in t( \tilde{G})$ and $ \tilde{g} \in t( \mbox{dom}(\tilde{f}))$
with $p \tilde{f} = f $ and $p\tilde{g} = g$, then $ \tilde{G}(fg) = \mbox{dom}(\tilde{f}
\tilde{g } ) = \mbox{dom}(\tilde{g})$, and $ (\tilde{G}f)g = (\mbox{dom}(\tilde{f}))
g = \mbox{dom}(\tilde{g})$, so $ \tilde{G}(fg) = (\tilde{G}f)g$. This function is an right
action.

Given any two objects $ \tilde{G}, \tilde{G'} \in \mbox{Ob}( \mathbf{H}_G)$, since
$\mathbf{\tilde{G}}$ is connected, there exists an arrow $ \tilde{f}: \tilde{G'}
\to \tilde{G}$, hence $\tilde{G'} = \mbox{dom}( \tilde{f}) = \tilde{G}(p\tilde{f})$.
This action is transitive.

(ii) If $f \in \pi( \mathbf{G}, G)$ such that $ \tilde{G}f = \tilde{G}$, then there exists
$\tilde{f} \in \pi( \mathbf{\tilde{G}}, \tilde{G})$ with $p\tilde{f} = f$, hence $f
\in p_\ast\pi( \mathbf{\tilde{G}}, \tilde{G})$. For the reverse direction, any $f
\in p_\ast\pi( \mathbf{\tilde{G}}, \tilde{G})$, there exists
$\tilde{f} \in \pi( \mathbf{\tilde{G}}, \tilde{G})$ with $p\tilde{f} = f$, since
 $t( \tilde{G})$ and $t(G)$ are bijective, this $ \tilde{f}$ is unique, hence $ \tilde{G}f
 = \mbox{dom}( \tilde{f}) = \tilde{G}$.

(iii) it is an easy result form group theory.
\end{proof}

If covering projection $p$ is $n$-fold, form Theorem~\ref {T: action} (iii) we know
$[ \pi( \mathbf{G}, G):
p_\ast\pi( \mathbf{\tilde{G}}, \tilde{G})] = n$.

A covering projection $p: \mathbf{\tilde{G}}  \to \mathbf{G}$ of groupoids is called
\emph{regular} if for every object $ \tilde{G}$ of $ \mathbf{\tilde{G}}$ the group
$p_\ast\pi( \mathbf{\tilde{G}}, \tilde{G})$ is a normal subgroup of
$\pi( \mathbf{G}, p\tilde{G})$. It is shown \cite [Theorem 9.6.4]{B} that when $p$ is a
connected covering projection $p$ is regular if and only if for given any loop in
$ \mathbf{G}$
(i.e., an element in any fundamental group of $ \mathbf{G}$), either every lifting or none
is loop in $ \mathbf{\tilde{G}} $.

If $ ( \mathbf{\tilde{G}}, p)$ is a covering groupoid of $ \mathbf{G}$, then a \emph{covering
transformation} is an isomorphism $h: \mathbf{\tilde{G}}  \to \mathbf{\tilde{G}}$ with
$ph = p$, that is, the following diagram commutes:
$$
\xymatrix{\mathbf{\tilde{G}}\ar[rr]^-h\ar[dr]_-p&&\mathbf{\tilde{G}}\ar[dl]^-p\\
&\mathbf{G}}
$$
Define Cov($\mathbf{\tilde{G}}/\mathbf{G}$) as the set of all covering transformations of
$ (\mathbf{\tilde{G}}, p)$. It's easy to see that Cov($\mathbf{\tilde{G}}/\mathbf{G}$) is a
group under composition of groupoid-morphisms.

\begin{theorem}\label {T: trans}
If $p: \mathbf{\tilde{G}} \to \mathbf{G}$ is a connected covering projection of groupoids,
then
$p$ is regular if and only if Cov($\mathbf{\tilde{G}}/\mathbf{G}$) acts transitively on
$\mbox{Ob}( \mathbf{H}_G)$ of the fibre $ \mathbf{H}_G$ over $G$.
\end{theorem}
\begin{proof} It's obvious that the group Cov($\mathbf{\tilde{G}}/\mathbf{G}$) acts on the
set $\mbox{Ob}( \mathbf{H}_G)$. We have to prove that when $p$ is regular the action is
transitive. Let $ \tilde{G}, \tilde{G'} \in \mbox{Ob}( \mathbf{H}_G)$, take any arrow
$a \in \mathbf{\tilde{G}}( \tilde{G}, \tilde{G'})$, then $a^{-1}\pi( \mathbf{\tilde{G}},
\tilde{G'})a = \pi(\mathbf{\tilde{G}}, \tilde{G})$. Since $p$ is regular,
$p_\ast(a^{-1}\pi( \mathbf{\tilde{G}}, \tilde{G'})a) =
p(a^{-1})p_\ast\pi( \mathbf{\tilde{G}}, \tilde{G'})pa =
p_\ast\pi( \mathbf{\tilde{G}}, \tilde{G'})
= p_\ast\pi(\mathbf{\tilde{G}}, \tilde{G})$. There exists a unique morphism
$h: \mathbf{\tilde{G}}  \to \mathbf{\tilde{G}}$ with
$ph = p$ and $h( \tilde{G'}) = \tilde{G}$ by the unique lifting theorem. Similarly,
there exists a unique morphism
$k: \mathbf{\tilde{G}}  \to \mathbf{\tilde{G}}$ with
$pk = p$ and $k( \tilde{G}) = \tilde{G'}$. By uniqueness $hk = 1$ and $kh = 1$, so $h$ is an
isomorphism and $h \in \mbox{Cov}(\mathbf{\tilde{G}}/\mathbf{G})$.

Conversely, assume that Cov($\mathbf{\tilde{G}}/\mathbf{G}$) acts transitively on
$\mbox{Ob}( \mathbf{H}_G)$: if $\tilde{G}, \tilde{G'} \in \mbox{Ob}( \mathbf{H}_G)$, there
exists $h \in \mbox{Cov}(\mathbf{\tilde{G}}/\mathbf{G})$ with $h( \tilde{G'}) = \tilde{G}$,
hence $h_\ast\pi( \mathbf{\tilde{G}}, \tilde{G'}) = \pi(\mathbf{\tilde{G}}, \tilde{G})$.
Since $p = ph$, then $p_\ast = p_\ast h_\ast$ and
$p_\ast\pi( \mathbf{\tilde{G}}, \tilde{G'}) =
p_\ast h_\ast\pi(\mathbf{\tilde{G}}, \tilde{G'})
= p_\ast\pi( \mathbf{\tilde{G}}, \tilde{G})$.
For any arrow $a \in \pi(\mathbf{G}, G)$, it has a lifting $ \tilde{a} \in t( \tilde{G'})$,
suppose $\mbox{dom}( \tilde{a}) = \tilde{G}$, i.e., $ \tilde{a} \in \mathbf{\tilde{G}}
(\tilde{G}, \tilde{G'})$, then $\tilde{a}^{-1}\pi( \mathbf{\tilde{G}},
\tilde{G'})\tilde{a} = \pi(\mathbf{\tilde{G}}, \tilde{G})$, hence
$$
\begin{aligned}
p_\ast(\tilde{a}^{-1}\pi( \mathbf{\tilde{G}}, \tilde{G'})\tilde{a})
&= p(\tilde{a}^{-1})p_\ast\pi( \mathbf{\tilde{G}}, \tilde{G'})p(\tilde{a})
= a^{-1}p_\ast\pi( \mathbf{\tilde{G}},\tilde{G'})a
= p_\ast\pi( \mathbf{\tilde{G}}, \tilde{G})\\
&= p_\ast\pi(\mathbf{\tilde{G}}, \tilde{G'})
\end{aligned} $$
so $p_\ast\pi(\mathbf{\tilde{G}}, \tilde{G'})$ is a normal group of $
\pi( \mathbf{G}, G)$, hence $p$ is regular.
\end{proof}

\begin{lemma}\label {L: free}
Let $ (\mathbf{\tilde{G}}, p)$ be a connected covering groupoid of $ \mathbf{G}$.
\begin{itemize}
\item [(i)] If $h \in \mbox{Cov}( \mathbf{\tilde{G}}/\mathbf{G})$ and $h \neq 1$, then $h$
has no fixed points.
\item [(ii)] If $h_1, h_2 \in \mbox{Cov}( \mathbf{\tilde{G}}/\mathbf{G})$ and there exists
$ \tilde{G} \in \mbox{Ob}( \mathbf{\tilde{G}})$ with $h_1( \tilde{G}) = h_2( \tilde{G})$,
then $h_1 = h_2$.
\end{itemize}
\end{lemma}
\begin{proof} These are easy results from the unique lifting theorem.
\end{proof}

\begin{corollary}
A connected covering projection $p$ is regular if and only if
$\mbox{Cov}( \mathbf{\tilde{G}}/\mathbf{G})$ acts principally on $\mbox{Ob}( \mathbf{H}_G)$
of the fibre $ \mathbf{H}_G$ over $G$. $\mbox{Cov}( \mathbf{\tilde{G}}/\mathbf{G})$ is
isomorphic to the quotient group $\pi( \mathbf{G}, G)/p_\ast\pi( \mathbf{\tilde{G}},
\tilde{G})$.
\end{corollary}
\begin{proof} The first statement follows from Theorem~\ref {T: trans} (transitivity) and
Lemma~\ref {L: free} (freeness). The second statement follows from Theorem~\ref {T: action} (ii)
and (iii).
\end{proof}

We can extend previous corollary to general case.

\begin{theorem}
Let $p: \mathbf{\tilde{G}} \to \mathbf{G}$ be a connected covering projection, $ \tilde{G}
\in \mbox{Ob}( \mathbf{H}_G)$ of the fibre $ \mathbf{H}_G$ over $G$, then
$$ \mbox{Cov}( \mathbf{\tilde{G}}/\mathbf{G}) \cong N_\pi(p_\ast\pi( \mathbf{\tilde{G}},
\tilde{G}))/p_\ast\pi( \mathbf{\tilde{G}},\tilde{G}).$$
where $N_\pi(p_\ast\pi( \mathbf{\tilde{G}}, \tilde{G}))$ is the normalizer subgroup of
$p_\ast\pi( \mathbf{\tilde{G}},\tilde{G})$ in the fundamental group $\pi( \mathbf{G}, G)$.
\end{theorem}
\begin{proof} At first, we construct a group homomorphism $\phi:
N_\pi(p_\ast\pi( \mathbf{\tilde{G}}, \tilde{G}))/p_\ast\pi( \mathbf{\tilde{G}},\tilde{G})
 \to \mbox{Cov}( \mathbf{\tilde{G}}/\mathbf{G})$.

Assume that any $[a] \in
N_\pi(p_\ast\pi( \mathbf{\tilde{G}}, \tilde{G}))/p_\ast\pi( \mathbf{\tilde{G}},\tilde{G})  $,
take a representative $a \in N_\pi(p_\ast\pi( \mathbf{\tilde{G}}, \tilde{G}))$, it has
a unique lifting $ \tilde{a} \in t( \tilde{G})$, denote its domain $\mbox{dom}( \tilde{a})$
by $ \tilde{G'}$, then $\pi( \mathbf{\tilde{G}}, \tilde{G'}) = \tilde{a}^{-1}
\pi( \mathbf{\tilde{G}}, \tilde{G})\tilde{a}$, hence
$p_\ast\pi( \mathbf{\tilde{G}}, \tilde{G'}) = p_\ast
\pi( \mathbf{\tilde{G}}, \tilde{G})$, since $N_\pi(p_\ast\pi( \mathbf{\tilde{G}}, \tilde{G}))$
is the normalizer subgroup. By unique lifting theorem, there exists a covering transformation
$h: \mathbf{\tilde{G}} \to \mathbf{\tilde{G}}$ with $h( \tilde{G'}) = \tilde{G}$. If $b\in [a]$,
then $b = ma $ where $m \in p_\ast\pi( \mathbf{\tilde{G}},\tilde{G})$, their liftings are
$\tilde{b} = \tilde{m}\tilde{a}$, so $\mbox{dom}( \tilde{b}) = \mbox{dom}( \tilde{a})$,
hence $b$ corresponds the same element that $a$ corresponds in the group
$\mbox{Cov}( \mathbf{\tilde{G}}/\mathbf{G})$. $ \phi$ is well-defined. It is obvious that
$\phi$ is a group homomorphism.

From the definition of $\phi$, if $\phi[b] = \phi[a]$, then $\mbox{dom}( \tilde{b}) =
\mbox{dom}( \tilde{a})$, where $ \tilde{b}$ and $ \tilde{a}$ are the unique liftings of $b$
and $a$ at $ \tilde{G}$ , respectively. Since $\mbox{codom}( \tilde{b}) =
\mbox{codom}( \tilde{a})
= \tilde{G}$, hence $ \tilde{b}\tilde{a}^{-1} \in \pi( \mathbf{\tilde{G}},\tilde{G})$,
and $ba^{-1} \in p_\ast\pi( \mathbf{\tilde{G}},\tilde{G}) $, then $[b] = [a]$.
$ \phi$ is injective.

Suppose that $g \in \mbox{Cov}( \mathbf{\tilde{G}}/\mathbf{G})$ and $g( \tilde{G'}) = \tilde{G}
$. Since $ \mathbf{\tilde{G}}$ is connected, there exists an arrow $ \tilde{a}:
\tilde{G'} \to \tilde{G}$. Hence $g_\ast\pi( \mathbf{\tilde{G}}, \tilde{G'})
= \pi( \mathbf{\tilde{G}}, g\tilde{G'}) = \pi(\mathbf{\tilde{G}}, \tilde{G})
= \tilde{a}\pi( \mathbf{\tilde{G}}, \tilde{G'}) \tilde{a}^{-1}$. Since $p = pg$, then
$p_\ast = p_\ast g_\ast$, and $p_\ast\pi( \mathbf{\tilde{G}}, \tilde{G'})
= p_\ast g_\ast\pi( \mathbf{\tilde{G}}, \tilde{G'})
= p_\ast\pi( \mathbf{\tilde{G}}, \tilde{G})
= (p\tilde{a})p_\ast\pi( \mathbf{\tilde{G}}, \tilde{G'}) (p\tilde{a}^{-1})
= (p\tilde{a})p_\ast\pi( \mathbf{\tilde{G}}, \tilde{G}) (p\tilde{a}^{-1})$.
Let $a = p\tilde{a}$, so $a \in N_\pi(p_\ast\pi( \mathbf{\tilde{G}}, \tilde{G}))$. Since
$\phi(a) \in \mbox{Cov}( \mathbf{\tilde{G}}/\mathbf{G})$ and $\phi(a)(\tilde{G'}) = \tilde{G}
$, by the uniqueness, $\phi(a) = g$. So $\phi$ is also surjective.
\end{proof}

\begin{corollary}
Let $( \mathbf{\tilde{G}}, p)$ be a universal covering groupoid of $ \mathbf{G}$. Then, for
$ G \in \mbox{Ob}( \mathbf{G})$, $\mbox{Cov}( \mathbf{\tilde{G}}/\mathbf{G})
 \cong \pi( \mathbf{G}, G)$.
\end{corollary}

Let $p: \mathbf{\tilde{G}} \to \mathbf{G} $ and $q: \mathbf{\tilde{H}} \to \mathbf{H}$ be
two connected universal covering projections. Let $f: \mathbf{H} \to \mathbf{G}$ be a morphism
of groupoids, a groupoid-morphism $\tilde{f}: \mathbf{\tilde{H}} \to \mathbf{\tilde{G}}$ is
said to be a \emph{covering of} $f$ if there is a commutative diagram
$$\xymatrix{
\mathbf{\tilde{H}}\ar[r]^-{\tilde{f}}\ar[d]_-q&\mathbf{\tilde{G}}\ar[d]^-p\\
\mathbf{H}\ar[r]_-f&\mathbf{G}}
$$

$f$ induces a group morphism $f_\ast: \pi(\mathbf{H}, H) \to \pi(\mathbf{G}, fH)$, define
$f_\#$ to be the bottom morphism of the following commutative diagram
$$\xymatrix{
\pi(\mathbf{H}, H)\ar[r]^-{f_\ast}\ar[d]_-\cong&\pi(\mathbf{G}, fH)\ar[d]^-\cong\\
\mbox{Cov}( \mathbf{\tilde{H}}/\mathbf{H})\ar[r]_-{f_\#}
&\mbox{Cov}( \mathbf{\tilde{G}}/\mathbf{G})}
$$

\begin{proposition}
Let $g \in \mbox{Cov}( \mathbf{\tilde{H}}/\mathbf{H})$ and let $\tilde{f}$ cover $f$, then
$f_\#(g)\tilde{f} = \tilde{f}g$. i.e., the following diagram commutes:
$$\xymatrix{
\mathbf{\tilde{H}}\ar[r]^-{g}\ar[d]_-{\tilde{f}}&\mathbf{\tilde{H}}\ar[d]^-{\tilde{f}}\\
\mathbf{\tilde{G}}\ar[r]_-{f_\#(g)}&\mathbf{\tilde{G}}}
$$
\end{proposition}
\begin{proof}
Both $f_\#(g)\tilde{f}$ and $ \tilde{f}g$ cover $f$. The equation holds by diagram chasing.
\end{proof}

Let $\Gamma$ be a group and $ \mathbf{G}$ a groupoid. An \emph{action} of $\Gamma$ on
$ \mathbf{G}$ assigns to each $\gamma \in \Gamma$ a morphism of groupoids $\gamma_\ast:
\mathbf{G} \to \mathbf{G}$ such that $1_\ast = 1: \mathbf{G} \to \mathbf{G}$, and if
$\gamma, \tau \in \Gamma$ then $(\gamma\tau)_\ast = \gamma_\ast \tau_\ast$. If $\gamma
\in \Gamma, G \in \mbox{Ob}( \mathbf{G})$ and $a \in \mbox{Mor}( \mathbf{G})$, write
$\gamma G$ and $\gamma a$ as the actions of $\gamma$ on $G$ and $a$, respectively.

An \emph{orbit groupoid} of the action is groupoid $ \mathbf{G}/\Gamma$ together with a
morphism $p: \mathbf{G} \to \mathbf{G}/\Gamma$ such that
\begin {itemize}
\item [(a)] If $ \gamma \in \Gamma$ and $g \in \mbox{Mor}( \mathbf{G})$, then
$ p(\gamma g) = p(g)$;
\item [(b)] The morphism $p$ is universal: if $f: \mathbf{G} \to \mathbf{H}$ is a morphism
of groupoids such that $f(\gamma g) = f(g)$ for all $ \gamma \in \Gamma$ and $g \in \mbox{Mor}
( \mathbf{G})$, then there is a unique morphism $f^\ast: \mathbf{G}/\Gamma \to \mathbf{H}$
of groupoids such that $f^\ast p = f$. In other words, the following triangle commutes:
$$\xymatrix{
\mathbf{G}\ar[rr]^p\ar[dr]_-f&& \mathbf{G}/\Gamma\ar[dl]^{\exists !f^\ast}\\
&\mathbf{H}}
$$
\end {itemize}
The morphism $p: \mathbf{G} \to \mathbf{G}/\Gamma$ is called an \emph{orbit morphism}.

$\mathbf{G}/\Gamma$ can be obtained from $ \mathbf{G}$ by imposing the relations $\gamma g
= g$ for all $\gamma \in \Gamma$ and all $g \in \mbox{Mor}( \mathbf{G})$. Brown gives a
construction of $\mathbf{G}/\Gamma$ in his book \cite [Setion 9.10.] {B}.

\begin{lemma}
If group $\Gamma$ acts freely on connected
groupoid $ \mathbf{G}$, by which we mean no non-identity
element of $ \Gamma$ fixes an object of $ \mathbf{G}$, then $p: \mathbf{G} \to
\mathbf{G}/\Gamma$ is a regular covering projection.
\end{lemma}
\begin{proof}
Group $\Gamma$ acts freely on groupoid $ \mathbf{G}$, this action induces two free actions of
$\Gamma$ on sets $\mbox{Ob}( \mathbf{G})$ and $\mbox{Mor}( \mathbf{G})$, respectively. Denote
the orbit of $G \in \mbox{Ob}( \mathbf{G})$ by $o(G) = \{\gamma G~ |~ \gamma \in \Gamma \}$
and the set of all orbits by $\mbox{Ob}( \mathbf{G})/\Gamma = \{ \mbox{ all orbits}~ o(G)\}$.
Let $ \mathbf{G}(o(G), o(G'))$ be the subset of $\mbox{Mor}( \mathbf{G})$ consisting of all
arrows $\{ a: C \to C' ~|~C \in o(G)~\mbox{and}~C'\in o(G')\}$, for any $\gamma \in \Gamma$,
$\gamma a: \gamma C \to \gamma C'$ is still in $ \mathbf{G}(o(G), o(G'))$, it is easy to see
that $ \Gamma$ also acts  freely on $ \mathbf{G}(o(G), o(G'))$. Denote
the orbit of $a \in \mathbf{G}(o(G), o(G'))$ by $o(a) = \{\gamma a~ |~ \gamma \in \Gamma \}$
and the set of all orbits by $\mathbf{G}(o(G), o(G'))/\Gamma = \{ \mbox{ all orbits}~ o(a)\}$.

By the define $ \mathbf{G}/\Gamma$ has as set $ \mbox{Ob}(\mathbf{G}/\Gamma)$
of objects the set
$\mbox{Ob}( \mathbf{G})/\Gamma$ and as set $\mathbf{G}/\Gamma(o(G), o(G'))$ of arrows from
$o(G)$ to $o(G')$ the set $\mathbf{G}(o(G), o(G'))/\Gamma$. The composite of arrows $o(a)
o(b) = o(\gamma(a)b)$ where $\gamma(\mbox{dom}(a)) = \mbox{codom}(b)$ and $\gamma$ is unique
since the action of $\Gamma$ is free. Suppose that $a' \in o(a)$ and $b' \in o(b)$, there exist
$\tau, \varsigma \in \Gamma$ such that $\tau a' = a$ and $\varsigma b' = b$, respectively.
similarly, $o(a')o(b') = o(\xi(a')b')$ where $\xi(\mbox{dom}(a')) = \mbox{codom}(b')$.
Since $\mbox{codom}(b) = \gamma(\mbox{dom}(a)) = \gamma\tau(\mbox{dom}(a'))
= \varsigma(\mbox{codom}(b')) = \varsigma\xi(\mbox{dom}(a'))$ and $\Gamma$ acts freely, then
$\gamma\tau = \varsigma\xi$. Hence $o(\gamma(a)b) = o(\gamma(\tau a')\varsigma b' )
= o(\varsigma\xi(a')\varsigma b') = o(\varsigma(\xi(a')b') = o(\xi(a')b')$, so the composite is
independent on the choice of representatives. The identity arrow is $o(1)$ and the inverse is
$o(a)^{-1} = o(a^{-1})$.

Define the orbit morphism $p: \mathbf{G} \to \mathbf{G}/\Gamma$ mapping the object $G \in
\mbox{Ob}( \mathbf{G})$ to $o(G)$ and arrow $a \in \mbox{Mor}( \mathbf{G})$ to $o(a)$. It is
obvious that $p$ is a morphism of groupoids. The map $p: t(G) \to t(o(G))$ is bijective: for
$a \in t(G)$ and $1\neq \gamma \in \Gamma$, $\gamma a \in t(\gamma G) $ where $t(G)
\cap t(\gamma G) = \emptyset$ since $\Gamma$ acts freely, so $p$ is injective. For every
$o(a) \in t(o(G))$, codomain of $a$ is in $o(G)$, hence there is $\gamma \in \Gamma$ such that
$\gamma(\mbox{codom}(a)) = G$, then $\gamma(a) \in t(G)$ with $o(\gamma(a)) = o(a)$, so
$p(\gamma(a)) = o(\gamma(a)) = o(a)$ and $p$ is surjective.
Hence $p$ is a covering projection.

Now $\Gamma \subset \mbox{Cov}(\mathbf{G}/(\mathbf{G}/\Gamma))$ because each
$\gamma \in \Gamma$ may be regarded as an isomorphism of $ \mathbf{G}$ with $p\gamma = p$.
Let $o(G) \in \mbox{Ob}( \mathbf{G}/\Gamma)$ and let $G \in \mbox{Ob}(\mathbf{G})$ be such that
$p(G) = o(G)$, as the objects of the fibre over $o(G)$ are $\{\gamma G: \gamma \in \Gamma \}$,
it follows that $\Gamma$, hence $\mbox{Cov}(\mathbf{G}/(\mathbf{G}/\Gamma))$, acts
transitively on these objects. By Theorem~\ref {T: trans}, $p$ is regular.
\end{proof}

\begin{lemma}
Let $( \mathbf{\tilde{G}}, p)$ be a regular connected
covering groupoid of $ \mathbf{G}$, and let
$ \Gamma = \mbox{Cov}( \mathbf{\tilde{G}}/\mathbf{G})$. There exists an isomorphism
$\varphi: \mathbf{G} \to \mathbf{\tilde{G}}/\Gamma$ making the following diagram commute:
$$\xymatrix{
\mathbf{\tilde{G}}\ar[d]_-p\ar[dr]^-q\\
\mathbf{G}\ar@{.>}[r]_-\varphi& \mathbf{\tilde{G}}/\Gamma}
$$
Moreover, $( \mathbf{\tilde{G}}, q)$ is a covering groupoid of $\mathbf{\tilde{G}}/\Gamma$.
\end{lemma}
\begin{proof}
Since $ \Gamma = \mbox{Cov}( \mathbf{\tilde{G}}/\mathbf{G})$, there is an action of
$ \Gamma$ on $ \mathbf{\tilde{G}} $ assigning to each $ \gamma \in \Gamma$ the corresponding
covering transformation. This action is free by Lemma~\ref {L: free}, hence
$( \mathbf{\tilde{G}}, q)$ is a regular covering groupoid of $\mathbf{\tilde{G}}/\Gamma$ by
previous lemma.

Since $p$ is a regular covering projection, $ \Gamma$ acts transitively on
$\mbox{Ob}( \mathbf{H}_G)$ of the fibre $ \mathbf{H}_G$ over $G$ by Theorem~\ref {T: trans}.
Take any $ \tilde{G} \in \mbox{Ob}( \mathbf{H}_G)$, hence $o( \tilde{G}) =
\mbox{Ob}( \mathbf{H}_G)$. Define $ \varphi(G) = o( \tilde{G})$,
it is independent on the choice
of $ \tilde{G}$, so it is well-define. Hence we can define $\varphi: \mbox{Ob}(\mathbf{G})
\to \mbox{Ob}(\mathbf{\tilde{G}}/\Gamma)$. It is easy to see that $\varphi$ is a bijection
on objects.

Let set of arrows $ \mathbf{\tilde{G}}(\mathbf{H}_G, \mathbf{H}_{G'}) = \{ \tilde{a}:
\tilde{G} \to \tilde{G'}~|~\tilde{G} \in \mbox{Ob}(\mathbf{H}_G)
~\mbox{and}~\tilde{G'}\in
\mbox{Ob}(\mathbf{H}_{G'})\}$, $\Gamma$ acts principally on
$ \mathbf{\tilde{G}}(\mathbf{H}_G, \mathbf{H}_{G'})$. It is easy to know that
$ \mathbf{\tilde{G}}(\mathbf{H}_G, \mathbf{H}_{G'})/\Gamma \cong \mathbf{G}(G, G')$.
Then $\varphi: \mathbf{G}(G, G') \to \mathbf{\tilde{G}}/\Gamma(o(\tilde{G}), o(\tilde{G')}$
defined by $a \mapsto o(\tilde{a})$, where $\tilde{a}$ is a lifting of $a$ in
$ \mathbf{\tilde{G}}(\mathbf{H}_G, \mathbf{H}_{G'})$, is a bijection on arrows.

It is obvious that $\varphi$ preserves identity and composite, so $\varphi$ is an isomorphism
of groupoids. It is also obvious that above triangle commutes.
\end{proof}

Let $(\mathbf{\tilde{G}}, p)$ and $( \mathbf{\tilde{H}}, q)$
be covering groupoids of $\mathbf{G}$
and $\mathbf{H}$, respectively. These covering groupoids are \emph{equivalent} if there exist
isomorphisms $\varphi$ and $\psi$ making the following diagram commute:
$$\xymatrix{
\mathbf{\tilde{H}}\ar[r]^-\varphi\ar[d]_-q&\mathbf{\tilde{G}}\ar[d]^p\\
\mathbf{H}\ar[r]_-\psi&\mathbf{G}}
$$
Hence in previous lemma the covering groupoids $(\mathbf{\tilde{G}}, p)$
and $(\mathbf{\tilde{G}}, q)$ are equivalent.

\begin{lemma}
Consider the commutative diagram of connected covering groupoids
$$\xymatrix{
\mathbf{\tilde{G}}\ar[dd]_-p\ar[dr]^-r\\
&\mathbf{\tilde{H}}\ar[dl]^-q\\
\mathbf{G}}
$$
where $p$ and $r$ are regular. Let $\Gamma = \mbox{Cov}( \mathbf{\tilde{G}}/\mathbf{G})$
and $\Pi = \mbox{Cov}( \mathbf{\tilde{G}}/\mathbf{\tilde{H}})$. Then there is a commutative
diagram
$$\xymatrix{
\mathbf{\tilde{G}}\ar[dd]_-{p'}\ar[dr]^-{r'}\\
&\mathbf{\tilde{G}}/\Pi\ar[dl]^-{q'}\\
\mathbf{\tilde{G}}/\Gamma}
$$
of covering groupoids, each of which is equivalent to the corresponding covering groupoids
in the original diagram.
\end{lemma}
\begin{proof}
By previous lemma, $(\mathbf{\tilde{G}}, r)$ is equivalent to $(\mathbf{\tilde{G}}, r')$
and $(\mathbf{\tilde{G}}, p)$ is equivalent to $(\mathbf{\tilde{G}}, p')$, where $r'$
and $p'$ are natural maps that send an object and an arrow in $\mathbf{\tilde{G}}$ into
their orbits, respectively.

If $\varphi: \mathbf{\tilde{G}} \to \mathbf{\tilde{G}}$ is an isomorphism with $r\varphi
= r$, then $p\varphi = qr\varphi = qr = p$, hence $\Pi =
\mbox{Cov}( \mathbf{\tilde{G}}/\mathbf{\tilde{H}}) \subset
 \Gamma = \mbox{Cov}( \mathbf{\tilde{G}}/\mathbf{G})$. For each object $\tilde{G} \in
 \mbox{Ob}(\mathbf{\tilde{G}})$, the $\Pi$-orbit of $\tilde{G}$ is contained in the
$\Gamma$-orbit of $\tilde{G}$. Define $q': \mathbf{\tilde{G}}/\Pi \to
\mathbf{\tilde{G}}/\Gamma$ to be the morphism that send a $\Pi$-orbit $o_\Pi(\tilde{G})$
of $\tilde{G}$ into the corresponding $\Gamma$-orbit $o_\Gamma(\tilde{G})$ of $\tilde{G}$
and an arrow $a \in \mathbf{\tilde{G}}/\Pi(o_\Pi(\tilde{G}),o_\Pi(\tilde{G}))$
to an arrow $q'a \in \mathbf{\tilde{G}}/\Gamma(o_\Gamma(\tilde{G}),o_\Gamma(\tilde{G}))$ such
that $a$ has a unique lifting $\tilde{a} \in t(\tilde{G'})$ and $q'a = p\tilde{a}$.
It is clear that $q'$ is a groupoid-morphism and $q'r' = p'$. Since $r'$ and $p'$ are
covering projections, $t(o_\Pi(\tilde{G})) \cong t(\tilde{G}) \cong t(o_\Gamma(\tilde{G}))$,
hence $q'$ is a covering projection.

In this diagram
$$
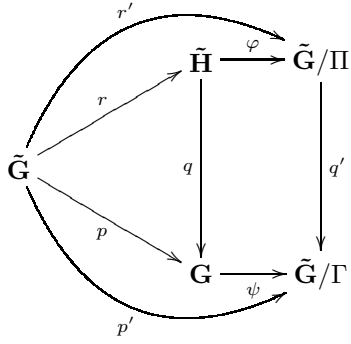
\xymatrix{
\\
&&\mathbf{\tilde{H}}\ar[r]^-\varphi\ar[dd]_-q& \mathbf{\tilde{G}}/\Pi\ar[dd]^-{q'}\\
\mathbf{\tilde{G}}\ar[urr]^-r\ar[drr]_-p\ar@/^3pc/[urrr]^-{r'}\ar@/_3pc/[drrr]_-{p'}\\
&&\mathbf{G}\ar[r]_-\psi&\mathbf{\tilde{G}}/\Gamma}
$$
$q'\varphi r = q'r' = p' = \psi p = \psi qr$. Since $r$ is a connected covering projection,
it is a groupoid-epi, hence $q'\varphi = \psi q$. From previous lemma, $\varphi$ and
$\psi$ are isomorphisms, so $(\mathbf{\tilde{H}}, q)$ is equivalent to
$(\mathbf{\tilde{G}}/\Pi, q')$.
\end{proof}

\begin{corollary}
Let $(\mathbf{\tilde{G}}, p)$ be a universal covering groupoid of $\mathbf{G}$, every covering
groupoid $(\mathbf{\tilde{H}}, q, r)$ of $\mathbf{G}$, where both $q: \mathbf{\tilde{H}}
\to \mathbf{G}$ and $r: \mathbf{\tilde{G}} \to \mathbf{\tilde{H}}$ are covering projections
with $qr = p$, is equivalent to
$(\mathbf{\tilde{G}}/\Gamma, \upsilon, p')$, where $\upsilon: \mathbf{\tilde{G}}/\Gamma
\to \mathbf{\tilde{G}}/\mbox{Cov}( \mathbf{\tilde{G}}/\mathbf{G})$
is a covering projection and $p': \mathbf{\tilde{G}} \to
\mathbf{\tilde{G}}/\Gamma$ is an orbit morphism
with $qr = p$, for some subgroup $\Gamma$ of
$\mbox{Cov}( \mathbf{\tilde{G}}/\mathbf{G})$.
\end{corollary}
\begin{proof}
There exists a morphism $r: \mathbf{\tilde{G}} \to \mathbf{\tilde{H}} $ of groupoids making
the following diagram commute:
$$
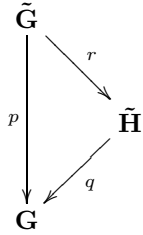
\xymatrix{
\mathbf{\tilde{G}}\ar[dd]_-p\ar[dr]^-r\\
&\mathbf{\tilde{H}}\ar[dl]^-q\\
\mathbf{G}}
$$
$(\mathbf{\tilde{G}}, r)$ is a covering groupoid of $\mathbf{\tilde{H}}$. Since
$\mathbf{\tilde{G}}$ is universal, both $p$ and $r$ are regular covering projections.
Therefore previous lemma implies that $(\mathbf{\tilde{H}}, q)$ is equivalent to
$(\mathbf{\tilde{G}}/\Gamma, \upsilon)$, where $\Gamma$ is a subgroup of
$\mbox{Cov}( \mathbf{\tilde{G}}/\mathbf{G})$.
\end{proof}

If $\Gamma$ is a group, let $\mbox{Sub}(\Gamma)$ be the family of all the subgroups of
$\Gamma$, and define a partial order $\Phi\preceq\Pi $ to mean $\Phi \subset \Pi$. Then
$\mbox{Sub}(\Gamma)$ is a lattice with $\Phi\vee\Pi$ the subgroup generated by $\Phi$ and
$\Pi$, and $\Phi\wedge\Pi = \Phi\cap\Pi$.

Let $(\mathbf{\tilde{G}}, p)$ be a universal covering groupoid of $\mathbf{G}$, and let
$\mbox{Lat}(\mathbf{\tilde{G}}/\mathbf{G})$ be the family of all equivalence classes of
covering groupoids $(\mathbf{\tilde{H}}, q, r)$
of $\mathbf{G}$. The previous corollary shows that each equivalence
class of covering groupoids can be expressed by representative
$(\mathbf{\tilde{G}}/\Gamma, \upsilon, p')$, where $\Gamma$ is a subgroup of
$\mbox{Cov}( \mathbf{\tilde{G}}/\mathbf{G})$. Define the partial order
$(\mathbf{\tilde{G}}/\Gamma, \upsilon, p') \preceq
(\mathbf{\tilde{G}}/\Pi, \omega, q')$ to mean there exists a covering
projection $s: \mathbf{\tilde{G}}/\Pi \to \mathbf{\tilde{G}}/\Gamma$ with $sq' = p'$. Then
$\mbox{Lat}(\mathbf{\tilde{G}}/\mathbf{G})$ is a lattice with $\mathbf{\tilde{G}}/\Gamma
\vee\mathbf{\tilde{G}}/\Pi$ the component of their pullback
$\mathbf{\tilde{G}}/\Gamma \times_\mathbf{G} \mathbf{\tilde{G}}/\Pi$ which $\mathbf{\tilde{G}}$
maps into, we still denote this component by
$\mathbf{\tilde{G}}/\Gamma \times_\mathbf{G} \mathbf{\tilde{G}}/\Pi$:
$$\xymatrix{
\mathbf{\tilde{G}}\ar[drr]^-{q'}\ar@{.>}[dr]|-{(p',q')}\ar[ddr]_-{p'}\\
&\mathbf{\tilde{G}}/\Gamma \times_\mathbf{G}\mathbf{\tilde{G}}/\Pi\ar[r]_-{\pi_2}\ar[d]
^-{\pi_1}& \mathbf{\tilde{G}}/\Pi\ar[d]^-{\omega}\\
&\mathbf{\tilde{G}}/\Gamma\ar[r]_-\upsilon & \mathbf{G}  }
$$
$\mathbf{\tilde{G}}/\Gamma
\wedge\mathbf{\tilde{G}}/\Pi$ their pushout $\mathbf{\tilde{G}}/\Gamma
\sqcup_\mathbf{\tilde{G}}
\mathbf{\tilde{G}}/\Pi$:
$$\xymatrix{\mathbf{\tilde{G}}\ar[r]^-{q'}\ar[d]_-{p'}
&\mathbf{\tilde{G}}/\Pi\ar[d]\ar[ddr]^-{\omega}
\\
\mathbf{\tilde{G}}/\Gamma \ar[r]\ar[drr]_-\upsilon& \mathbf{\tilde{G}}/\Gamma
\sqcup_\mathbf{\tilde{G}}
\mathbf{\tilde{G}}/\Pi\ar@{.>}[dr]|-{(\upsilon, \omega)}\\
&&\mathbf{G}}$$

\begin{theorem}\label {T: main}
Let $(\mathbf{\tilde{G}}, p)$ be a universal covering groupoid of $\mathbf{G}$ and
$\Gamma = \mbox{Cov}( \mathbf{\tilde{G}}/\mathbf{G})$, then
\begin {itemize}
\item [(i)] The function $\gamma: \mbox{Lat}(\mathbf{\tilde{G}}/\mathbf{G})
\to \mbox{Sub}(\Gamma)$, defined by $(\mathbf{\tilde{H}}, q, r) \mapsto
\mbox{Cov}( \mathbf{\tilde{G}}/\mathbf{\tilde{H}})$, is an order reversing bijection with
inverse $\delta: \Pi \to (\mathbf{\tilde{G}}/\Pi, \upsilon, p')$.
\item [(ii)] $(\mathbf{\tilde{G}}/\mbox{Cov}( \mathbf{\tilde{G}}/\mathbf{\tilde{H}}), q)
= (\mathbf{\tilde{H}}, q)$ and $ \mbox{Cov}( \mathbf{\tilde{G}}/(
\mathbf{\tilde{G}}/\Pi)) = \Pi$.
\item [(iii)] $$\mathbf{\tilde{G}}/\Pi\vee\Phi = \mathbf{\tilde{G}}/\Pi
\sqcup_\mathbf{\tilde{G}}\mathbf{\tilde{G}}/\Phi,$$
$$\mathbf{\tilde{G}}/\Pi\wedge\Phi = \mathbf{\tilde{G}}/\Pi
\times_\mathbf{{G}}\mathbf{\tilde{G}}/\Phi,$$
$$\mbox{Cov}( \mathbf{\tilde{G}}/\mathbf{\tilde{H}}\vee\mathbf{\tilde{K}})
= \mbox{Cov}( \mathbf{\tilde{G}}/\mathbf{\tilde{H}})\cap
\mbox{Cov}( \mathbf{\tilde{G}}/\mathbf{\tilde{K}}),$$
$$\mbox{Cov}( \mathbf{\tilde{G}}/\mathbf{\tilde{H}}\wedge\mathbf{\tilde{K}})
= \mbox{Cov}( \mathbf{\tilde{G}}/\mathbf{\tilde{H}})\vee
\mbox{Cov}( \mathbf{\tilde{G}}/\mathbf{\tilde{K}}).$$
\item [(iv)] Fold of $\mathbf{\tilde{H}} = [\Gamma:
\mbox{Cov}( \mathbf{\tilde{G}}/\mathbf{\tilde{H}})]$ and $[\Gamma: \Pi] = $ fold of
$\mathbf{\tilde{G}}/\Pi$.
\item [(v)] $(\mathbf{\tilde{H}}, p)$ is a regular covering groupoid of $\mathbf{G}$ if and
only if $\mbox{Cov}( \mathbf{\tilde{G}}/\mathbf{\tilde{H}})$ is a normal subgroup of
$\Gamma$.
\end {itemize}
\end{theorem}
\begin{proof}
(i) If $\Pi \subset \Gamma = \mbox{Cov}(\mathbf{\tilde{G}}/\mathbf{G})$, then $\gamma\delta
(\Pi) = \mbox{Cov}(\mathbf{\tilde{G}}/(\mathbf{\tilde{G}}/\Pi))$, call this last group
$\Pi^\ast$. Then $\Pi^\ast$ consists of all isomorphisms $h: \mathbf{\tilde{G}} \to
\mathbf{\tilde{G}}$ making the following diagram commute:
$$\xymatrix{
\mathbf{\tilde{G}}\ar[rr]^-h\ar[dr]_-\upsilon&&\mathbf{\tilde{G}}\ar[dl]^-\upsilon\\
&\mathbf{\tilde{G}}/\Pi}
$$
where $\upsilon: \mathbf{\tilde{G}} \to \mathbf{\tilde{G}}/\Pi$ is the natural covering
projection. If $\theta \in \Pi$ and $\tilde{G}, \tilde{G'} \in \mbox{Ob}(\mathbf{\tilde{G}}),
\tilde{a} \in \mathbf{\tilde{G}}(\tilde{G}, \tilde{G'})$, then $\upsilon(\tilde{G}) =
\upsilon(\theta\tilde{G}), \upsilon(\tilde{a}) = \upsilon(\theta\tilde{a})$ by the definition
of $\mathbf{\tilde{G}}/\Pi$, hence $\upsilon\theta = \upsilon$, so $\theta \in \Pi^\ast$
 and $\Pi \subset \Pi^\ast$. For the reverse inclusion, let $h \in \Pi^\ast$ with $\upsilon h
 = \upsilon$. If $\tilde{G} \in \mbox{Ob}(\mathbf{\tilde{G}})$ then the orbit $o(\tilde{G} )
 = o(h(\tilde{G}))$. By the definition of $\Pi$-orbit, there exists $\theta \in \Pi$ with
 $\theta(h(\tilde{G})) = \tilde{G}$. Since $\theta \in \Pi \subset \Pi^\ast$, it follows
 that $\theta h \in \Pi^\ast$. By Lemma~\ref {L: free} (ii), $\theta h = 1$, and $h = \theta
 ^{-1} \in \Pi$, hence $\Pi^\ast \subset \Pi$. So $\gamma\delta(\Pi) = \Pi$.

$\delta\gamma$ is the composite $(\mathbf{\tilde{G}}/\Pi, \upsilon) \mapsto
\mbox{Cov}(\mathbf{\tilde{G}}/(\mathbf{\tilde{G}}/\Pi)) = \Pi^\ast \mapsto (\mathbf{\tilde{G}}
/\Pi^\ast, \upsilon^\ast)$. Since there exists $\Pi = \Pi^\ast$, hence $\upsilon = \upsilon^
\ast$, $\delta\gamma$ is identity. So $\gamma$ (or $\delta$) is a bijection.

When $\mathbf{\tilde{H}} \preceq \mathbf{\tilde{K}}$ in
$\mbox{Lat}(\mathbf{\tilde{G}}/\mathbf{G})$, there exists morphism $r: \mathbf{\tilde{K}}
\to \mathbf{\tilde{H}}$ of groupoids making the following diagram commute:
$$\xymatrix{
\mathbf{\tilde{G}}\ar[dd]_-t\ar[dr]^-s\\
&\mathbf{\tilde{K}}\ar[dl]^-r\\
\mathbf{\tilde{H}}}
$$
hence $\mbox{Cov}( \mathbf{\tilde{G}}/\mathbf{\tilde{K}}) \subset
\mbox{Cov}( \mathbf{\tilde{G}}/\mathbf{\tilde{H}})$, $\gamma$ is order reversing.

(ii) This is just the statement that $\delta\gamma$ and $\gamma\delta$ are identity functions.

(iii) These equations follows from the basic properties of order-reversing bijection of
lattices and the fact that both $\gamma$ and $\delta = \gamma^{-1}$ are order-reversing
bijection.

(iv) Since $(\mathbf{\tilde{G}}, p)$ is a universal covering of $\mathbf{G}$, $\Gamma
= \mbox{Cov}( \mathbf{\tilde{G}}/\mathbf{G}) \cong \pi(\mathbf{G}, G)$, $
\mbox{Cov}( \mathbf{\tilde{G}}/\mathbf{\tilde{H}}) \cong \pi(\mathbf{\tilde{H}}, \tilde{H}),
\Pi = \mbox{Cov}( \mathbf{\tilde{G}}/(\mathbf{\tilde{G}}/\Pi)) \cong \pi(
(\mathbf{\tilde{G}}/\Pi), \ast)$. This statement follows from the formula of fold (Theorem~
\ref {T: action} (iii)).

(v) This is just the definition of regular covering.
\end{proof}

\begin{remark}
We have a classification of covering groupoids over a fixed groupoid $\mathbf{G}$ in terms of
subgroups of covering transformation group $\mbox{Cov}( \mathbf{\tilde{G}}/\mathbf{G })$.
Higgins \cite [Proposition 36.] {H}
gives a classification of covering groupoids over $\mathbf{G}$ in terms of
subgroups of fundamental group $\pi(\mathbf{G}, G)$. Since
$\mbox{Cov}( \mathbf{\tilde{G}}/\mathbf{G }) \cong \pi(\mathbf{G}, G)$, the two classification
results coincide.
\end{remark}

\begin{corollary}
Let $(\mathbf{\tilde{G}}, p)$ be a universal covering groupoid of $\mathbf{G}$. If $\Pi$ is
a subgroup of $\mbox{Cov}( \mathbf{\tilde{G}}/\mathbf{G }) (\cong \pi(\mathbf{G}, G))$, then
$$\pi(\mathbf{\tilde{G}}/ \Pi, \ast) \cong \Pi.$$
\end{corollary}
\begin{proof}
From the proof (iv) of previous theorem there exist isomorphisms
$\pi(\mathbf{\tilde{G}}/ \Pi, \ast) \cong
\mbox{Cov}( \mathbf{\tilde{G}}/(\mathbf{\tilde{G} }/\Pi)) = \Pi$.
\end{proof}


\begin{thebibliography}{4}

\bibitem{B} R. Brown, \emph{Topology: A geometric accout of general
topology, homotopy types and the fundamental groupoid}, John Wiley
\& sons, 1988.


\bibitem{H}  P.J. Higgins, \emph{Notes on
Categories and Groupoids}, Van Nostrand Reinhold Company, 1971.


\bibitem{J} P.T. Johnstone, \emph{Topos theory}, Academic press, 1977.


\bibitem{M-M} S. Mac Lane and I. Moerdijk, \emph{Sheaves in geometry and logic},
Springer-Verlag, 1992.

\bibitem{R} J.J. Rotman, \emph{An introduction to algebraic topology}, Graduate Texts in
Mathematics 119, Springer-Verlag, 1988.


\end{thebibliography}
\end{document}